\documentclass[runningheads]{llncs}
\usepackage{amsmath}
\usepackage{amsfonts}
\usepackage{amssymb}
\usepackage{graphicx}

\usepackage{multirow}
\usepackage{wrapfig}
\usepackage{floatflt} 
\usepackage{float}

\usepackage{tikzpagenodes}
\usepackage{tikz}
\usetikzlibrary{arrows.meta}
\usetikzlibrary{arrows,matrix,positioning}

\usepackage{pgfplots}
\pgfplotsset{compat=1.17}
\usepgfplotslibrary{groupplots}

\begin{document}
	\title{Towards Efficient Time Stepping for \\Numerical Shape Correspondence}
	\titlerunning{Towards Efficient Time Stepping for Numerical Shape Correspondence}
	%
	\author{Alexander Köhler\inst{1}
	        \and Michael Breuß\inst{1}
		}
	\authorrunning{Alexander Köhler and Michael Breuß}
	%
	\institute{Institute for Mathematics, Brandenburg Technical University, 
	           Platz der Deutschen Einheit 1, 03046 Cottbus, Germany\\
		\email{\{koehlale,breuss\}@b-tu.de}}
	\maketitle              
	\begin{abstract}
    The computation of correspondences between shapes is a principal task in shape analysis.
    To this end, methods based on partial differential equations (PDEs) have been established, 
    encompassing e.g.\ the classic heat kernel signature as well as numerical solution schemes 
    for geometric PDEs. In this work we focus on the latter approach. 
    
    We consider here several time stepping schemes. The goal of this investigation 
    is to assess, if one may identify a useful property of methods for time integration for the shape analysis
    context. Thereby we investigate the dependence on time step size, since the
    class of implicit schemes that are useful candidates in this context should ideally
    yield an invariant behaviour with respect to this parameter.
    
    To this end we study integration of heat and wave equation on a manifold. In order to facilitate
    this study, we propose an efficient, unified model order reduction framework for these models.
    We show that specific $l_0$ stable schemes are favourable for numerical shape analysis. 
    We give an experimental evaluation of the methods at hand of classical TOSCA data sets.
    \end{abstract}
	%
	%
	%
	
	\section{Introduction}
	The computation of shape correspondences is a fundamental task in computer vision with many 
	potential applications, cf. \cite{Kaick2011,Bronstein2008a}. In the setting of three-dimensional 
	shape analysis, the underlying problem amounts to identify an explicit relation between the 
	surface elements of two or more shapes. The variety of possible shape correspondence mappings that
	is of interest in applications includes non-rigid transformations where shapes are just almost 
	isometric, allowing e.g.\ to match different poses of human or animal shapes.
	
	An important solution strategy is to achieve a pointwise shape correspondence using so called 
	descriptor based methods. For this, a feature descriptor has to be computed that characterizes 
	each point on a shape by describing the surrounding shape surface geometry. A mathematically 
	sound approach to compute such shape signatures is to make use of the spectral decomposition 
	of the Laplace-Beltrami operator, see e.g. \cite{Rustamov2007,Sun2009,Bronstein2008a}.
	To this end, the Laplace-Beltrami operator may be incorporated in a certain variety of partial differential
	equations (PDEs) that are potentially useful as models for feature computation.
	The arguably most important classic signature is the heat kernel signature (HKS) \cite{Sun2009}
	which relies on the heat equation, however also versions of the
	Schrödinger equation leading to the wave kernel signature (WKS) \cite{Aubry2011} as well 
	as the hyperbolic wave equation \cite{Dachsel2017} have already been proposed.
	
	For the computation of such PDE-based signatures, the basic task amounts to resolve the underlying PDEs
	on a manifold representing a shapes' boundary. The HKS and WKS both rely on the eigenfunction expansion
	of the Laplace-Beltrami operator to tackle this task and for achieving efficient algorithms.
	An alternative to the spectral approach is to consider the numerical integration of the underlying PDEs 
	as proposed in \cite{Dachsel2017,Dachsel2017a,Baehr2018}. However, when following this path it is highly advocated 
	to employ very efficient computational means such as the model order reduction framework presented in 
	\cite{Baehr2018,Baehr2019a} in order to avoid high computational times. 
	
	In previous works based on numerical integration, first order implicit time integration has been studied in detail \cite{Dachsel2017a}. 
	For the wave equation model, it has been shown that backward differencing in time may yield favorable results over simple central differences \cite{Dachsel2017}. Especially, as has been illustrated in \cite{Dachsel2017}, the classic wave equation may yield in some experiments 
	results of higher correspondence quality compared to the other mentioned models. Let us note that this holds again in particular when assessing the models using numerical integration with implicit first order time stepping. 
	
	{\bf Our contributions.}
	In this paper we give an account of ongoing work on improving numerical PDE-based shape descriptors with respect
	to computational efficiency. The underlying question addressed in this paper is, if
	the class of implicit schemes has particularly useful key properties that should be respected
	in the shape analysis context. We present a study of three implicit time stepping schemes for 
	heat and wave equation on manifolds, respectively. We show that standard finite differences may not be appropriate since the typical initial condition used for the construction of shape signatures, which is a discrete Dirac delta function, may yield oscillatory artefacts that may spoil correspondence quality. In order to resolve this issue, we propose to adopt $l_0$-stable methods. For performing the study, we have developed a unified 
	numerical model order reduction framework based on \cite{Baehr2018}. We give an experimental account of our investigations at hand of selected shape data sets.

	\section{Modelling of the Shape Correspondence Framework}
	
	The shape of three-dimensional geometric object can be described by its bounding surface, given as
	a compact two-dimensional Riemannian manifold $\mathcal{M}\subset \bbbr^3$. 
	Two such shapes $\mathcal{M}$ and $\widetilde{\mathcal{M}}$ are considered to be isometric if there is a smooth homeomorphism between the corresponding object surfaces that preserves the intrinsic distances between surface points.
	For many applications, isometry may only hold approximately, leading to the notion of almost isometric
	shapes.
	
	As indicated, a widely used descriptor class that can handle almost isometric transformations relies on  PDEs on manifolds. We will consider in this paper the corresponding heat and wave equation, respectively. 
	These PDEs may be put together in a more general form as described below.
	
	All the mentioned PDEs rely on the Laplace-Beltrami operator as a fundamental building block.
	This operator is the geometric version of the Laplace operator, meaning that it takes into account the local 
	curvature of a smooth manifold in 3D. Applied to a scalar-valued function $u$, it can be written as: 
	\begin{equation}
		\Delta_{\mathcal{M}} u = \frac{1}{\sqrt{\vert g \vert}} \sum_{i,j=1}^{2} \partial_i \left( \sqrt{\vert g \vert}g^{ij} \partial_j u\right)
	\end{equation}
	where $\vert g \vert$ is the determinant of the metric tensor $g\in \bbbr^{2\times 2}$ that describes locally the geometry, and $g^{ij}$ are the entries of its inverse, see e.g.\ \cite{DoCarmo2016} for more details on the differential geometric notions.
	
	\paragraph{Heat Equation} 
	The heat equation on a manifold is a proper shape descriptor~\cite{Sun2009} and reads as
	\begin{equation}\label{eq:heat_raw}
	\partial_t u(x,t) = \Delta_{\mathcal{M}} u(x,t).\quad x \in \mathcal{M},~t \in I
	\end{equation}
	This PDE basically describes heat flow along a surface $\mathcal{M}$.
	 
	\paragraph{Wave Equation} Proposed in~\cite{Dachsel2017} for computing shape correspondences, 
	this equation can be written as  
	\begin{equation} \label{eq:wave_raw}
		\partial_{tt} u(x,t) = \Delta_{\mathcal{M}} u(x,t).\quad x \in \mathcal{M},~t \in I
	\end{equation}
	
	\paragraph{Meta PDE}
	We can combine the PDEs from above into a single equation as
	\begin{equation}
		\label{eq:core}
		\phi \partial_{tt} u(x,t) + \psi \partial_t u(x,t) = \Delta_{\mathcal{M}} u(x,t) .\quad x \in \mathcal{M},~t \in I
	\end{equation}
	We see if we choose $\phi = 0$ and $\psi=1$ we receive heat equation (\ref{eq:heat_raw}) and for $\phi = 1$ and $\psi=0$ it will lead us to the wave equation (\ref{eq:wave_raw}). The PDE formulation \eqref{eq:core} is useful for us in this paper, since it allows us to describe in a compact way numerical developments in later sections.
	
	\paragraph{Initial condition for feature computation}
	Both equations need initial conditions to be solved accurately. The classic initial condition used for PDE-based feature generation is the Dirac Delta function $u(x,0)=u_0(x)= u_{x_i}$ centred around a point $x_i \in \mathcal{M}$, cf.\ \cite{Sun2009}. We also adopt this initial condition in the numerical integration framework. Because equation (\ref{eq:wave_raw}) is a PDE of second order in time we need an additional condition at the initial velocity. A canonical choice is to consider the zero initial velocity condition $ \partial_t u(x,0) = 0$. 
	
	
	
	\paragraph{Feature Descriptor}
	For establishing shape correspondence we consider for each point on an objects' surface a feature descriptor, which is a computational object that contains geometric shape information. Within our framework, a simple way to construct it is to restrict the spatial component of $u(x,t)$ to
	\begin{equation}
		f_{x_i}(t):=u(x,t)\vert_{x=x_i} \quad \mbox{with} \quad u(x,0)\vert_{x=x_i} = u_{x_i}
		\label{eq:feature_descriptor}
	\end{equation}
	and call $f_{x_i}$ the feature descriptor at the location $x_i \in \mathcal{M}$. 
	
	With this definition we may adopt a physical interpretation of \eqref{eq:feature_descriptor}. 
	The feature descriptor $f_{x_{i}}$ describes heat transferred away from $x_{i}$ 
	respectively the motion amplitudes of a wave front emitted and observed again at $x_{i}$. 
	Thereby the processes run on the surface $\mathcal{M}$, as indicated via the Laplace-Beltrami operator. 
	
	
	\paragraph{Shape Correspondence}
	To compare the feature descriptors for different locations $x_i \in \mathcal{M}$ and $\tilde{x}_j \in \widetilde{\mathcal{M}}$ on two different shapes $\mathcal{M}$ and $\widetilde{\mathcal{M}}$, we define a distance $d_f(x_i, \tilde{x}_j)$ via the $L_1$ norm 
	\begin{equation}
		d_f(x_i, \tilde{x}_j) = \int\limits_{I} \vert f_{x_i} - f_{\tilde{x}_j} \vert \textmd{d}t 
		\label{eq:dist_int}
	\end{equation}
	The tuple of locations $ (x_i, \tilde{x}_j) \in \mathcal{M} \times \widetilde{\mathcal{M}}$  with the smallest distance belong together.
	
	\paragraph{Basic Discretisation}

\begin{floatingfigure}[tr]{0.42\textwidth}
   \centering
		\begin{tikzpicture}[scale=0.75]
			\pgftext{\includegraphics[width=0.425\textwidth]{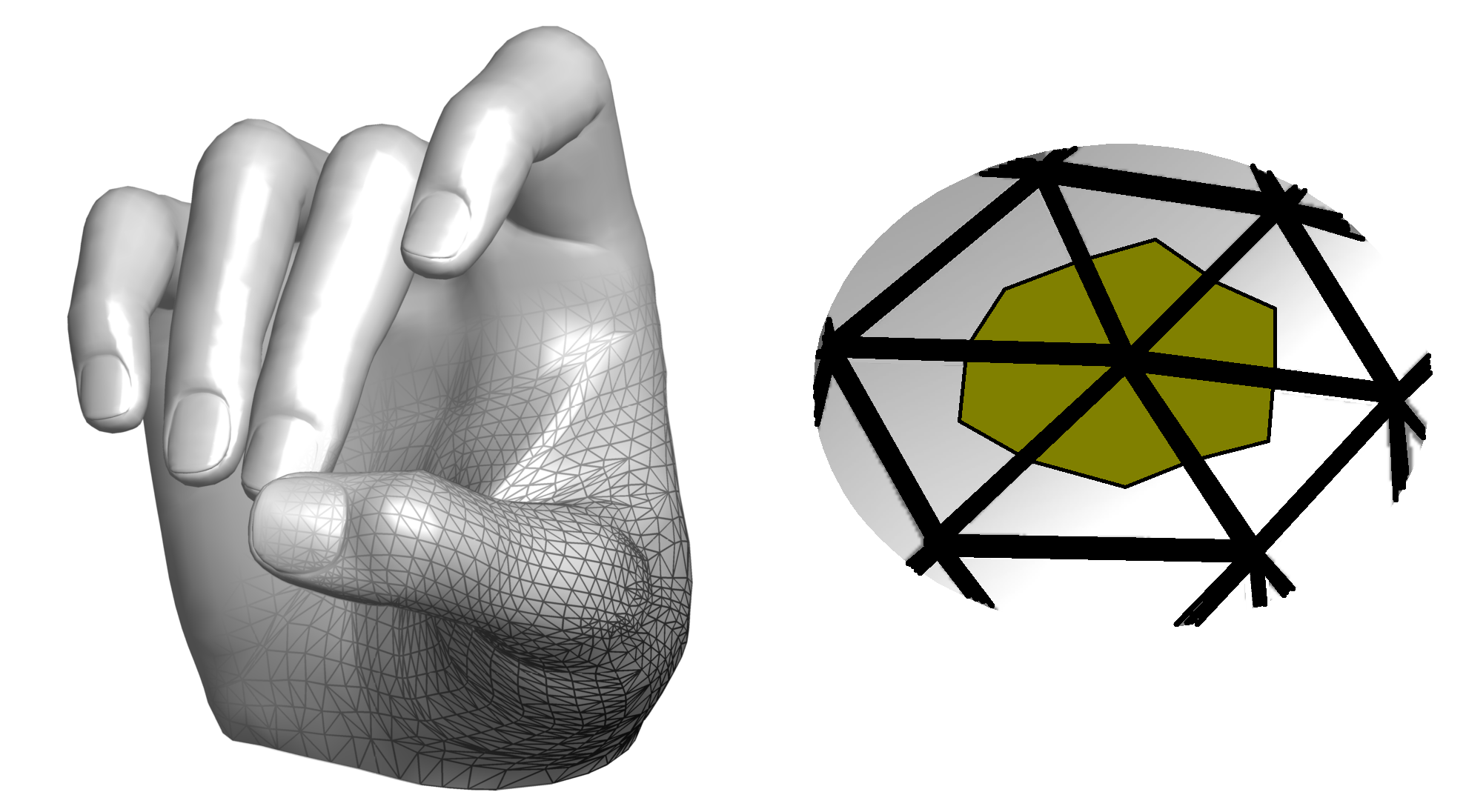} }  at (0pt,0pt); 
			\draw[ultra thick] (1.34,0.06) ellipse (1.1cm and 0.85cm);
			\draw[ultra thick] (-0.5,0) -- (0.55,0.65);
			\draw[ultra thick] (-0.5,-0.2) -- (0.8,-0.682);
		\end{tikzpicture}
		\caption{Continuous $\&$ discrete shape\label{fig:mesh}}
\end{floatingfigure}

	The discrete surface representation is given by a set of points $P := (x_1, \ldots, x_N)$ containing coordinates of each point $x_i$ and one for triangles $T$. $T$ stores triples of points forming a triangle. We denote this data set tuple as $\mathcal{M}_d = (P,T)$. 
	Furthermore we denote by $\Omega_i$ the barycentric cell volume surround the point $x_i$.
	The discrete shape is given by non-uniform linear triangles.
	Compare also Figure \ref{fig:mesh} (figure adopted from \cite{Baehr2018}).

	\section{Spatial Discretisation}
	We employ the equation (\ref{eq:core}) with a finite volume method for discretisation. 
	First, we consider (\ref{eq:core}) over a cell $\Omega_i$ and a time interval $I_k$ so that the integration in time and space will lead to 
	\begin{equation}
		\int \limits_{I_k}  \int \limits_{\Omega_i}   \phi \partial_{tt} u(x,t) + \psi \partial_t u(x,t) \, \mathrm{d}x \,  \mathrm{d}t  
		=  \int \limits_{I_k}  \int \limits_{\Omega_i}  \Delta_\mathcal{M} \,  u(x,t)   \, \mathrm{d}x \, \mathrm{d}t \label{he1}
	\end{equation}
	The definition of the cell average 
	\begin{equation}
		u_i(t)=u(\bar{x}_i,t)= \frac{1}{|\Omega_i|} \int \limits_{\Omega_i}  u(x,t)  \,   \mathrm{d}x
	\end{equation}
	where $|\Omega_i|$ denotes the area of $\Omega_i$, is used to define the averaged Laplacian as
	\begin{equation}
		L u_i(t)= \frac{1}{|\Omega_i|} \int \limits_{\Omega_i}   \Delta_\mathcal{M} \,u(x,t) \, \mathrm{d}x  
	\end{equation}
	Evaluation of this integral making use of the divergence theorem will lead to a line integral over the boundary of each 
	cell. Application of the cotangent weight scheme \cite{Meyer2003} turns the resulting integrals 
	over a full shape into the ODE system
	\begin{equation}
		\phi \mathbf{\ddot{u}}(t) + \psi \mathbf{\dot{u}}(t) = L\mathbf{u}(t) \label{eq:ode_all}
	\end{equation} 

	The cotangent weight scheme will give us the discrete Laplace-Beltrami operator $L\in \bbbr^{N \times N}$ via the sparse matrix representation $L = D^{-1} W$. The matrix $W$ represents spatial weights computed
	from the cotangent formulae, encoding spatial connectivity and local geometry information.
	The matrix $D=\text{diag}(|\Omega_1|,\ldots,|\Omega_i|,\ldots,|\Omega_N|) $ contains the local cell areas and all functions $u_i(t)$ are stored in 
	the $N$-dimensional vector $\mathbf{u}(t) = (u_1(t), \ldots, u_N(t))^{\top}$.
	
	Revisiting equation (\ref{eq:ode_all}), we may notice that for $\phi=0$ we have an ODE system of first order. 
	Additionally we want to transform the equation for $\phi \neq 0$ into an ODE system of first order as well, since this 
	is advantageous for the proceeding. We find
	\begin{equation}
		\mathbf{\dot{q}}(t) = H\mathbf{q}(t)
		\label{eq:wave_system}
	\end{equation}
	with $\mathbf{q}(t) = \big(\mathbf{u}(t),\mathbf{\dot{u}}(t)\big)^{\top} \in  \bbbr^{2N}$,  
	\begin{equation}
		H = \begin{pmatrix}
			0 & I \\ \frac{1}{\phi} L & -\frac{\psi}{\phi}I
		\end{pmatrix} \stackrel{\phi=1}{=} \begin{pmatrix}
			0 & I \\ L & -\psi I
		\end{pmatrix} \in \bbbr^{2N \times 2N}
	\end{equation}
	and $I \in \bbbr^{N\times N}$ the identity matrix. 
	
	Let us turn to the discrete initial conditions. Choosing  $\mathbf{\dot{u}}(0) = \mathbf{0}$ let us only need to consider the initial spatial conditions $\mathbf{u}(0)$.  Again using cell averages,
%
	the normalised initial condition at the location $x_i$ can be written as 
	\begin{equation}
		\mathbf{u_{x_i}}=\mathbf{u}(x_i,0)=(0,\dots,0,|\Omega_i|^{-1},0,\dots, 0)^\top \label{eq:initial_cond_for_u}
	\end{equation} 

	\paragraph{Eigenproblem and Modal Coordinate Reduction}
		Systems like (\ref{eq:ode_all}) have to deal with large sparse matrices and consequently high cost of computation time. With the Model Order Reduction (MOR) we can approximate our high dimensional system with lower dimensional one. In \cite{Nouri2014,Qu2004} one may find an overview about the MOR topic. 
	    For this work we are interested in the specific Modal Coordinate Reduction (MCR) technique
	    as presented in \cite{Baehr2019a}, which we briefly recall now.  
	
	Solving equation (\ref{eq:ode_all}) via MCR, it is essential to calculate the eigenvalues from the 
	Laplace-Beltrami Operator $L \mathbf{v} = \lambda \mathbf{v}$.
    With $L = D^{-1}W$ we can transform the eigenvalue problem to the generalised eigenvalue Problem
	\begin{equation}
		W \mathbf{v} = \lambda D \mathbf{v} \label{eq:general_eigenproblem}
	\end{equation}  
	By the properties of $W$ and $D$ all eigenvalues are real and the eigenvectors are linear independent. 
	Additionally they are $D$-orthogonal with $\mathbf{v}_i^{\top} D \mathbf{v}_j = \delta_{ij}$. Thus we 
	obtain the equalities
	\begin{equation}
		I = V^{\top} D V, \qquad L = V \Lambda V^{\top} D, \qquad \Lambda = V^{\top} W V
	\end{equation}
	with $V$ being the right eigenvector matrix of $L$ and $\Lambda$ the diagonal matrix of eigenvalues.   
 	
 	To execute the MCR transformation we take advantage of these considerations and substitute equation (\ref{eq:ode_all}) with $\mathbf{u} = V \mathbf{w}$ to obtain:
 	\begin{equation}
 		\phi V\ddot{\mathbf{w}} + \psi V\dot{\mathbf{w}}  = LV\mathbf{w}
 	\end{equation}
 	Multiplication with $V^{\top} D$ from the left results in 
 	\begin{equation}
 		\phi\ddot{\mathbf{w}} + \psi\dot{\mathbf{w}}  = \Lambda \mathbf{w} \label{eq:almost_reduced}
 	\end{equation}
 		Equation (\ref{eq:almost_reduced}) is still the full system. To reduce it we are only interested in the first $r \ll N$ ordered eigenvalues $0 = \lambda_1 < \lambda_2 \leq \ldots \leq \lambda_r$.
 		This gives us the reduced matrix $\Lambda_r \in \bbbr^{r \times r}$ extracted from $\Lambda$ and $V_r \in \bbbr^{N \times r}$. Consequently we obtain the reduced model 
 	\begin{equation}
 		\phi\ddot{\mathbf{w}}_r + \psi \dot{\mathbf{w}}_r  = \Lambda_r \mathbf{w}_r \quad \mbox{where} \quad \mathbf{w}_r = V_r^{\top} D \mathbf{u}  \label{eq:reduced_model}
 	\end{equation} 
 	Analogously to the procedure by which we obtained equation (\ref{eq:wave_system}) from (\ref{eq:ode_all}), we can build 
 	\begin{equation}
 		\mathbf{\dot{p}}_r(t) = H_r\mathbf{p}_r(t) \quad \mbox{with} \quad H_r =
 		\begin{pmatrix}
 			0 & I_r \\ \Lambda_r & -\psi I_r
 		\end{pmatrix} \in \bbbr^{2r \times 2r} \label{eq:reduced_system}
  	\end{equation}
 	and $\mathbf{p}_r = (\mathbf{w}_r, \dot{\mathbf{w}}_r)^{\top} \in  \bbbr^{2r}$.
 	
 	Now, we want to discuss the transformation of the initial conditions. Remember that $ \mathbf{u}(x_i,0) = (0,\dots,0,|\Omega_i|^{-1},0,\dots, 0)^\top$  and together with the second part of (\ref{eq:reduced_model}) we obtain $\mathbf{w}_r(x_i,0) = V_r^{\top} \mathbf{e}_i$. Nothing will change for the initial velocity $\mathbf{\dot{u}}(0) = \mathbf{0}$, because we set it zero and so will $\mathbf{\dot{w}}_r(0) = \mathbf{0}$.
 	

	\section{Time Discretisation}

	For discretisation of the explored time interval we divide $\left[0,t_M\right]$ into intervals $I_k = \left[t_{k-1},t_k\right]$ with $\tau = t_k-t_{k-1}$ as time increment, and we set $t_0 = 0$.  
	
\begin{figure}[t]
	    \centering
	    \includegraphics[width=\textwidth]{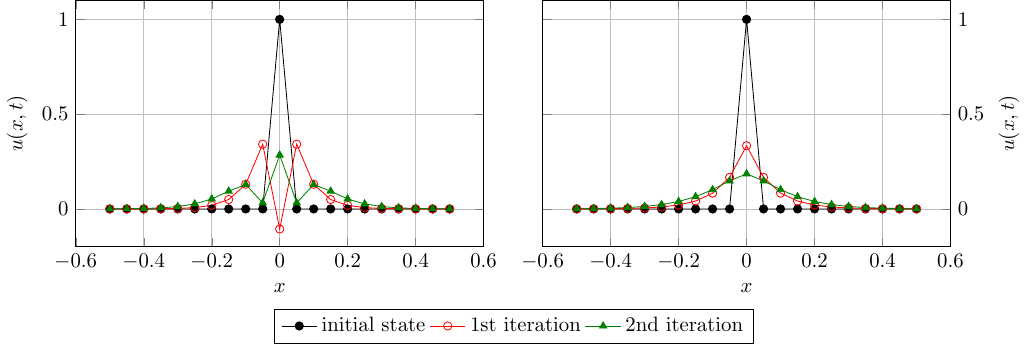}
	    \caption{This first two iterations of a heat equation solved with the Crank-Nicolson method (left) and with the $l_0$-stable implicit Euler scheme (right), respectively.}
	    \label{fig:cn_vs_ie}
	\end{figure}

	Let us now turn to time derivatives.
	The classic method of second order in time accuracy is the Crank-Nicolson method \cite{Morton2005}.
	We conjecture at this point that this method may not be beneficial in our shape analysis framework. 
	Illustrating this, we calculate the first two iterations of a heat equation with the ($l_0$ unstable) Crank-Nicolson method in comparison to the $l_0$-stable implicit Euler method, cf. Figure \ref{fig:cn_vs_ie}. 
	For the solution of the heat equation we expect a smooth propagation away from the peak. 
	However, the Crank-Nicolson method causes significant oscillations in the solution, even negative
	values appear. 
	We conjecture at this point that such oscillations may spoil correspondence quality, when
	performing numerical shape correspondence. This example illustrates why we aim to explore 
	$l_0$ stability of schemes.  
	
	
	\paragraph{Second Order $l_0$-stable Scheme}
	In this section, we briefly explain the idea and approach of the method of \cite{Twizell1996}. We have chosen this method because it does not depend on complicated Butcher tableaus, but offers an intuitive approach.
	
	Looking at equation (\ref{eq:wave_system}) or (\ref{eq:reduced_system}), we find these are of the form 
	\begin{equation}
		\dot{\mathbf{x}}(t) = A \mathbf{x}(t)
	\end{equation}
	The solution of this is $\mathbf{x}(t) = \exp(tA)$ or for a time step further 
	\begin{equation}
		\mathbf{x}(t+\tau) = \exp((t+\tau)A) = \exp(\tau A) \mathbf{x}(t)
	\end{equation}  
	Calculating $\exp(\tau A)$ is not a trivial task, so we use an approximation 
	\begin{equation}
		R(\tau A) \approx \exp(\tau A)
	\end{equation}
	instead. With the approximation 
	$R(\tau H) = (I-a\tau H)^{-1}(I+(1-a)\tau H)$
	we obtain a first-order method. Choosing $a$ will lead to different well-known schemes. For $a=0$ we receive the explicit method and for $a=1$ we obtain the implicit method, like it is used in \cite{Baehr2018,Baehr2019a} for time discretisation. The Crank-Nicolson scheme can be obtained by setting $a = 0.5$. 
	
	However, two out of three are not of second order and, in case of the explicit and Crank-Nicolson method, are not $l_0$-stable. Twizell and co-authors \cite{Twizell1996} present with
	\begin{equation}
		R(\tau A) = (I-r_1 \tau A)^{-1}(I-r_2 \tau A)^{-1} (I+(1-a)\tau A)
		\label{second-order-stable}
	\end{equation}
	with $r_{1,2} = \frac{1}{2}\left( a \mp \sqrt{a^2-4a + 2} \right)$ 
	a method of second order and choosing $a = 2-\sqrt{2}-\varepsilon$, with $\varepsilon$ an arbitrarily small positive number, will provide an $l_0$-stable method. 
	
	These considerations lead to the following discrete procedure for the wave equation:
	\begin{equation}
		\mathbf{p}^{k+1} = (I_{2r}-r_1 \tau H_r)^{-1}(I_{2r}-r_2 \tau H_r)^{-1} (I_{2r}+(1-a)\tau H_r) \mathbf{p}^k \label{eq:wave_eq_all_l0_scheme}
	\end{equation}
	with $\mathbf{p}^k = \mathbf{p}(t_k)$ and $\mathbf{p}^{k+1} = \mathbf{p}(t_{k+1}) = \mathbf{p}(t+\tau)$. 
	Analogously for the heat equation, we receive 
	\begin{equation}
		\mathbf{w}^{k+1} = (I_r-r_1 \tau \Lambda_r)^{-1}(I_r-r_2 \tau \Lambda_r)^{-1} (I_r+(1-a)\tau \Lambda_r) \mathbf{w}^k \label{eq:heat_eq_l0_scheme}
	\end{equation}
	As discussed in \cite{Baehr2019a} we adopt the temporal domain $\left[0, t_M\right]$ for the heat equation to $\left[0, t^*_{\text{h}}\right]$ and for wave and damped wave equation to $\left[0, t^*_{\text{w}}\right]$ by
	\begin{equation}
	    t^*_{\text{h}}(\lambda_r) = \frac{t_M \sqrt{\lvert \lambda_N \rvert}}{\sqrt{ \lvert \lambda_r \rvert }} 
	    \qquad \text{resp.} \qquad
	    t^*_{\text{w}}(\lambda_r) = \frac{t_M \sqrt[4]{\lvert \lambda_N \rvert}}{\sqrt[4]{ \lvert \lambda_r \rvert }} \label{eq:t_star}
	\end{equation}
	The computational parameters are thus $\tau = \frac{t^*}{M}$, with $M$ being the number of iterations.  

	\section{Experiments}
	
	\paragraph{Hit Rate and Geodesic Error} 
	The percentage Hit Rate is defined as $TP/(TP+FP)$, where $TP$ is the number of true positives and $FP$ is the number of false positives events. 
	
	The Princeton benchmark protocol \cite{Kim2011} delivers the important key points for the evaluation of the correspondence quality. It defines the quality of a matching by the normalised intrinsic distance $d_{\mathcal{M}}(x_i,x^*)/\sqrt{A_\mathcal{M}}$ from the calculated matching $x_i$ to the ground-truth correspondence $x^*$. A matching is accepted to be true if the normalised intrinsic distance is smaller then the threshold $0.25$.  
	
	\paragraph{Dataset}
	For the experiments we used a selction of shapes from the TOSCA dataset \cite{Bronstein2008a} containing a total of 80 shapes of animals and humans in different poses. The range of shape resolution in TOSCA
	goes from 4\,344 points for the wolf shape to 52\,565 points for the david and michael shape. 
	

	
	\paragraph{Experimental Settings}
	
	    \begin{figure}[t]
        \centering
        \includegraphics{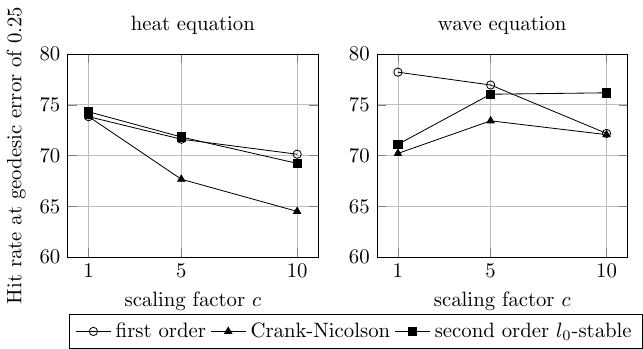}
		
		
        \caption{\label{fig:time_step_scaling}
        The figure displays mean hit rates in \% at geodesic error of 0.25 
        for different shapes and schemes in dependence of the scaling factor $c$.
        The shape sets over which we average are cat, centaur, dog, horse, and wolf, taken from the TOSCA
        Dataset \cite{Bronstein2008a}.
        }
    \end{figure}

	As discussed, we use three methods to solve the heat and wave equation on different shapes. 
	As a reference method we employ the $l_0$-stable implicit Euler method of first order in time \cite{Baehr2019a}. 
	We compare it with the discussed second order methods, i.e. Crank-Nicolson 
	and the $l_0$-stable method \eqref{second-order-stable} from \cite{Twizell1996}. 
	
	Let us recall that we aim to explore the possibility to increase the time step size $\tau$.
	The point with this is, that an implicit scheme should theoretically perform stable independently
	on the time step size \cite{Morton2005}. Thus we expect that a rescaling of time step size may reveal
	if the $l_0$ stability is consistently a useful property.
	
	The time step $\tau$ and the number of iterations depend directly on each other through $\tau M = t^* $. 
	Since we fixed the value for $t^*$ with the calculations from (\ref{eq:t_star}) an increase of $\tau$ leads to a decrease of $M$. We denote the factor by which we have increased the time step size $\tau$ by $c$. We have chosen $c = \left\lbrace 1,5,10 \right\rbrace$ for our experiments. With these values for $c$ the number of iterations changed to $M = \left\lbrace 100,20,10 \right\rbrace$.

	\paragraph{Results for the Heat Equation}
	A plot of results can be found in the left figure within Figure \ref{fig:time_step_scaling}. 
	We see that all three methods lose some accuracy when the time step size is increased. 
	The Crank-Nicolson method gives the worst hit rates. For $c=1$ it can still keep up with the other methods (73.84\%), but it loses significantly more accuracy (64.56\%) when the time step size is increased, 
	compared to the other two methods. The second order $l_0$-stable method and the implicit Euler method are very similar. 
	Generally speaking, if we compare the hit rates from the $l_0$-stable methods (first and second order), 
	we observe nearly no differences. 

\begin{figure}[ht] 
        \centering
        \includegraphics[scale = 0.9]{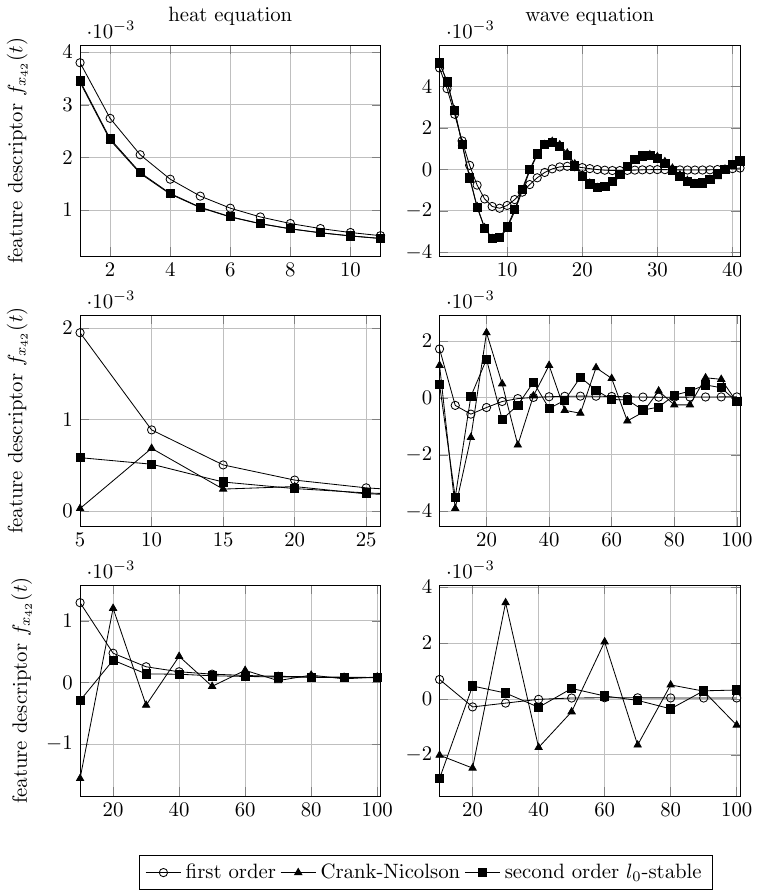}
        \caption{\label{fig:fd_compare_schemes}
        The feature descriptors of the wolf shape in dependencies of the time. The feature descriptor is received from the heat (left) or wave (right) equation solved with a first order method \cite{Baehr2019a}, the Crank-Nicolson method and a $l_0$-stable second order method \cite{Twizell1996}. From Top to Bottom we see the different time scaling factor $c$. Starting with $c=1$ (top) going to $c= 5$ (middle) and end up with $c=10$ (bottom). In the top plots the second-order method are close to each other so they appear as one curve. 
        }
    \end{figure}

	
	A look at the left plots of Figure \ref{fig:fd_compare_schemes} shows that the feature descriptor for the individual methods is not very different. However, the beneficial aspect of an $l_0$-stable method 
	becomes evident, since for an increase of the time step the Crank-Nicolson method produces 
	clear oscillations, like in the simple motivation example in Figure \ref{fig:cn_vs_ie}. 
	
	The results for the heat equation emphasize the usefulness of the $l_0$ stability as an important notion of practical importance in numerical shape analysis.

 	\paragraph{Results for Wave Equation}
    The right figure within Figure \ref{fig:time_step_scaling} is providing the plots of the mean hit 
    rates. 
    Interesting is the fact, that with choosing a second order method we lose about $8\%$ accuracy in the hit rate for $c=1$ compared to the Euler integrator. At this point, further investigation is needed.
    However, we can notice an increase in hit rate again by increasing the time step size. 
    Again we observe here that $l_0$ stability is a beneficial property.
    
     Looking at the right plots in Figure \ref{fig:fd_compare_schemes}, we see that a much more detailed feature descriptor can be produced by the second order methods. However, this does not increase the hit rate. 
    
    Summarising the findings, we observe again that the $l_0$ stable schemes perform better than the
    classical Crank-Nicolson method.

	\section{Conclusion and Further Work}
	
	Summarising the experiments, we have validated that $l_0$ stability is an important property 
	for numerical integration in shape analysis. This is surely a useful aspect for possible
	future developments.
	
	There are some other indirect findings, e.g. that the use of higher order schemes for time integration 
	may not be of highly beneficial impact, since the $l_0$-stable second order time integrator did not 
	perform significantly better than the implicit Euler scheme. However, this may be partly due to the 
	fact of the observed drop in accuracy in case of the wave equation, a subject we aim to resolve 
	in future work.
	
	\bibliographystyle{splncs04}
	\bibliography{bibliography}

\begin{thebibliography}{10}
\providecommand{\url}[1]{\texttt{#1}}
\providecommand{\urlprefix}{URL }
\providecommand{\doi}[1]{https://doi.org/#1}

\bibitem{Aubry2011}
Aubry, M., Schlickewei, U., Cremers, D.: The wave kernel signature: {A} quantum
  mechanical approach to shape analysis. In: 2011 {IEEE} International
  Conference on Computer Vision Workshops ({ICCV} Workshops). {IEEE} (Nov 2011)

\bibitem{Baehr2018}
B{\"{a}}hr, M., Breu{\ss}, M., Dachsel, R.: {Fast Solvers for Solving Shape
  Matching by Time Integration} (2018)

\bibitem{Baehr2019a}
B{\"{a}}hr, M., Breu{\ss}, M., Dachsel, R.: {Model Order Reduction for Shape
  Correspondence} arXiv: 1910.00914v2 (2020)

\bibitem{Bronstein2008a}
Bronstein, A.M., Bronstein, M.M., Kimmel, R.: Numerical {G}eometry of
  {N}on-{R}igid {S}hapes. Springer-Verlag GmbH (2008)

\bibitem{Dachsel2017a}
Dachsel, R., Breu{\ss}, M., Hoeltgen, L.: {Shape Matching by Time Integration
  of Partial Differential Equations}. In: Lecture Notes in Computer Science,
  pp. 669--680. Springer International Publishing (2017)

\bibitem{Dachsel2017}
Dachsel, R., Breu{\ss}, M., Hoeltgen, L.: {The Classic Wave Equation Can Do
  Shape Correspondence}. In: Computer Analysis of Images and Patterns, pp.
  264--275. Springer International Publishing (2017)

\bibitem{DoCarmo2016}
Do~Carmo, M.P.: Differential geometry of curves and surfaces: revised and
  updated second edition. Dover Publications (2016)

\bibitem{Kaick2011}
van Kaick, O., Zhang, H., Hamarneh, G., Cohen-Or, D.: {A Survey on Shape
  Correspondence}. Computer Graphics Forum  \textbf{30}(6),  1681--1707 (Jul
  2011)

\bibitem{Kim2011}
Kim, V.G., Lipman, Y., Funkhouser, T.: Blended intrinsic maps. In: {ACM}
  {SIGGRAPH} 2011 papers on - {SIGGRAPH} '11. {ACM} Press (2011)

\bibitem{Meyer2003}
Meyer, M., Desbrun, M., Schr{\"{o}}der, P., Barr, A.H.: Discrete
  {D}ifferential-{G}eometry {O}perators for {T}riangulated 2-{M}anifolds. In:
  Mathematics and Visualization, pp. 35--57. Springer Berlin Heidelberg (2003)

\bibitem{Morton2005}
Morton, K.W., Mayers, D.F.: {Numerical Solution of Partial Differential
  Equations}. Cambridge University Press (Apr 2005)

\bibitem{Nouri2014}
Nouri, S.: {Advanced Model-Order Reduction Techniques for Large Scale Dynamical
  Systems}. Ph.D. thesis, Ottawa-Carleton (2014)

\bibitem{Qu2004}
Qu, Z.Q.: {Model Order Reduction Techniques with Applications in Finite Element
  Analysis}. Springer London (2004)

\bibitem{Rustamov2007}
Rustamov, R.M.: {Laplace-Beltrami Eigenfunctions for Deformation Invariant
  Shape Representation, functional map}. In: Proceedings of the Fifth
  Eurographics Symposium on Geometry Processing. pp. 225--233. SGP
  {\textquoteright}07, Eurographics Association, Goslar, DEU (2007)

\bibitem{Sun2009}
Sun, J., Ovsjanikov, M., Guibas, L.: {A Concise and Provably Informative
  Multi-Scale Signature Based on Heat Diffusion}. Computer Graphics Forum
  \textbf{28}(5),  1383--1392 (Jul 2009)

\bibitem{Twizell1996}
Twizell, E.H., Gumel, A.B., Arigu, M.A.: {Second-order,{$L_0$}-stable methods
  for the heat equation with time-dependent boundary conditions}. Advances in
  Computational Mathematics  \textbf{6}(1),  333--352 (Dec 1996)

\end{thebibliography}
\end{document}